\documentclass[conference]{IEEEtran}

\usepackage{cite}
\usepackage{amsmath,amssymb,amsfonts}
\usepackage{algorithmic}
\usepackage{graphicx}
\usepackage{textcomp}
\usepackage{xcolor}

\usepackage{float}

\usepackage{listings}
\usepackage[hyphens]{url} 

\usepackage{lipsum}
\usepackage{enumitem}
\usepackage{comment}

\setlist[description]{
  labelindent=0cm,
  leftmargin=1.8cm, 
  labelwidth=1.5cm,
  style=standard,
  itemsep=2pt,      
  parsep=0pt
}

\usepackage{hyperref}

\begin{document}
\title{Electricity Price-Aware Scheduling of Data Center Cooling\\
\thanks{This material is based upon work supported by the U.S. Department of Energy’s Office of Energy Efficiency and Renewable Energy (EERE) under the Wind Energy Technologies Office (WETO) Award Number DE-EE0011269, the Massachusetts Clean Energy Center and the Maryland Energy Administration. The views expressed herein do not necessarily represent the views of the U.S. Department of Energy, the United States Government, the Massachusetts Clean Energy Center or the Maryland Energy Administration.}
}

\author{
\IEEEauthorblockN{Arash Khojaste, Jonathan Pearce, Golbon Zakeri}
\IEEEauthorblockA{\textit{Department of Mechanical and Industrial Engineering}\\
\textit{University of Massachusetts Amherst}\\
Amherst, MA, USA\\
akhojaste@umass.edu, jonathanpear@umass.edu, gzakeri@umass.edu}

\and
\IEEEauthorblockN{Yuanrui Sang}
\IEEEauthorblockA{\textit{Department of Electrical and Computer Engineering} \\
\textit{University of Massachusetts Amherst}\\
Amherst, MA, USA\\
ysang@umass.edu}
}

\maketitle

\begin{abstract}
Data centers are becoming a major consumer of electricity on the grid, with cooling accounting for about 40\% of that energy. As electricity prices vary throughout the day and year, there is a need for cooling strategies that adapt to these fluctuations to reduce data center cooling costs. In this paper, we present a model for electricity price-aware cooling scheduling using a Markov Decision Process (MDP) framework to reliably estimate the cooling system's operational costs and facilitate investment-phase decision-making. We utilize Quantile Fourier Regression (QFR) fits to classify electricity prices into different regimes while capturing both daily and seasonal patterns. We simulate 14 years of operation using historical real-time electricity price and outdoor temperature data, and compare our model against heuristic baselines. The results demonstrate that our approach consistently achieves lower cooling costs. This model is useful for grid operators interested in demand response programs and data center investors looking to make investment decisions.
\end{abstract}

\vspace{0.4em}


\section*{Nomenclature}
\vspace{0.4em}
\textbf{Sets and Indices}
\begin{description}
    \item[$t \in \mathcal{T}$] Discrete time index; $\mathcal{T} = \{1, \dots, N\}$.
    \item[$\theta \in \Theta$] Discrete indoor air temperature (°C); $\Theta = \{\theta_{\min}$, $\dots, \theta_{\max}\}$.
    \item[$p \in \mathcal{P}$] Electricity price regime index; $\mathcal{P} = \{1, \dots, M\}$.
    \item[$a \in \mathcal{A}$] Number of active chillers; $\mathcal{A} = \{0, \dots, A_{\max}\}$.
\end{description}
\vspace{0.4em}

\textbf{Parameters}
\begin{description}

    \item[$N$] Total number of discrete time steps.
    \item[$M$] Total number of electricity price regimes.
    \item[$A_{\max}$] Maximum number of active chillers.    
    
    \item[$\theta_{\min}$] Minimum indoor temperature used to define $\Theta$.
    \item[$\theta_{\max}$] Maximum indoor temperature used to define $\Theta$.
    
    \item[$T_{\min}$] Lower bound of indoor temperature safety range (°C).
    \item[$T_{\max}$] Upper bound of indoor temperature safety range (°C).
    
    \item[$v^{\text{under}}_{t}$] Temperature violation below $T_{\min}$ at time $t$ (°C).
    \item[$v^{\text{over}}_{t}$] Temperature violation above $T_{\max}$ at time $t$ (°C).

    \item[$\gamma_{\text{env}}$] Heat transfer coefficient for ambient exchange (W/°C).
    \item[$\Delta t$] Duration of a time step, e.g., one hour (s).
    \item[$\rho_{\text{air}}$] Air density (kg/m$^3$).
    \item[$c_{p,\text{air}}$] Specific heat of air (J/kg·°C).
    \item[$V_{\text{room}}$] Volume of the data center (m$^3$).
    \item[$C_{\text{air}}$] Thermal capacitance of indoor air (J/°C).
    \item[$A_{\text{floor}}$] Floor area of the data center (m$^2$).
    \item[$h_{\text{slab}}$] Slab thickness (m).
    \item[$\rho_{\text{concrete}}$] Concrete density (kg/m$^3$).
    \item[$c_{p,\text{concrete}}$] Specific heat of concrete (J/kg·°C).
    \item[$C_{\text{slab}}$] Thermal capacitance of slab (J/°C).
    \item[$C_{\text{equipment}}$] Total estimated thermal capacitance of IT equipment and infrastructure (J/°C).
    \item[$C_{\text{heat}}$] Total effective thermal capacitance (J/°C).
    \item[$Q^{\text{base}}$] Total base heat load on idle from IT equipment (W).
    \item[$i_t$] Number of active CPU cores at time $t$.
    \item[$\phi$] Heat produced per active core (W/core).
    \item[$Q_t$] Total heat load at time $t$ (W).
    \item[$T^{\text{out}}_t$] Outdoor temperature at time $t$ (°C).
    \item[$\eta$] Cooling rate per unit chiller level (W/level).
    \item[$\text{COP}_t$] Coefficient of performance of the cooling system at time $t$.
    \item[$e_{t,a}$] Electricity consumed at time $t$ under action $a$ (kWh).
    \item[$\pi_{t,p}$] Electricity price at time $t$ in regime $p$ (\$/MWh).
    \item[$\mathbb{P}^{t}_{p \rightarrow p'}$] Transition probability between electricity price regimes.
    
    \item[$\lambda_{\text{under}}$] Penalty per °C deviation below $T_{\min}$ (\$/°C).
    \item[$\lambda_{\text{over}}$] Penalty per °C deviation above $T_{\max}$ (\$/°C).

    \item [$c_{t,s,a}$] Immediate cost of being in state $s$ at time $t$ and taking action $a$ (\$).
    \item [$P_{t, s, s', a}$] Transition probability from $s = (\theta, p)$ to $s' = (\theta', p')$ under action $a$ at time $t$.
\end{description}

\vspace{0.4em}

\textbf{Decision Variables}
\begin{description}
    \item[$x_{t,s,a}$] Probability of being in state $s = (\theta, p)$ and taking action $a$ at time $t$.
\end{description}


\section{Introduction}
\vspace{-0.2em}

The data center industry is expanding rapidly, making it one of the fastest-growing industries worldwide. In the U.S., demand for data centers continues to rise due to increased computational needs driven by artificial intelligence, cloud computing, and big data. This growth significantly increases data centers' energy consumption \cite{NREL2024cooling}.
Based on estimates, data centers are expected to consume about 9\% of U.S. annual electricity generation by 2030, with cooling systems alone accounting for up to 40\% of that consumption \cite{EPRI2023AI}.
Given the substantial share of cooling in total energy consumption, optimizing data center cooling systems can significantly reduce energy costs and usage \cite{zhang2021survey}. Having an accurate estimate of optimized cooling costs well before the construction stage of a data center is important for both engineering design and financial planning. In this paper, we address how to estimate these costs, accounting for the stochastic nature of electricity prices and the variability of seasonal and daily temperatures. We utilize pre-cooling strategies that can be used to respond flexibly to electricity price fluctuations.

Recently, there have been improvements in cooling architectures, such as hybrid liquid-air systems and in-row cooling, which improve thermal efficiency. However, they are mostly hardware-driven optimizations. While Zhou et al. \cite{zhou2024advancements} provide an overview of these trends, Zhang et al. \cite{zhang2021survey} provide a comprehensive survey on different strategies used for data center cooling systems.
They categorize cooling strategies by their sources into \textbf{free cooling} and \textbf{mechanical refrigeration}. Free cooling uses ambient conditions (e.g., through air-side or water-side economizers, heat pipes), while mechanical refrigeration relies on active components such as chillers or compressors.
In the delivery phase, mechanical refrigeration may deliver cooling to the IT load via air (e.g., CRAC/CRAH) or liquid (e.g., cold plates or chilled water).

Jin et al. \cite{jin2024climate} analyze optimal cooling configurations across various climate zones, showing that free cooling requires low outdoor temperatures. Therefore, it is not effective during the summer and during that time, data centers must rely on mechanical refrigeration. This makes mechanical refrigeration essential and indispensable for the warmer seasons. \cite{zhou2024advancements}
In the model setup section of this paper, we assume free cooling is sufficient during cooler months and focus on mechanical refrigeration during the summer. This setup aligns with many data centers nowadays in the Massachusetts area as well (see \cite{MGHPCC}). However, our model is general enough to support year-round planning as well.

On optimizing the operational planning of the data center cooling systems, there are two main approaches popular in the literature: Model Predictive Control (MPC) (see \cite{lazic2018datacenter}, \cite{parolini2011cyber}) and Reinforcement Learning (RL) (see \cite{li2015multigrid}, \cite{zhang2018practical}, \cite{mnih2015human}, \cite{le2021deep}, \cite{sutton1998introduction}). While each has its advantages and challenges (see \cite{zhang2021survey}), they are specifically designed for real-time decision-making rather than for investment-phase decisions. In contrast, our model provides an estimate of the operational costs and operational plan associated with cooling a data center based on its physical specifications. It is particularly designed for the investment phase, helping two key stakeholders: grid operators and data center investors. For grid operators, it helps them assess the potential of demand response programs and works as a method to estimate the value of demand response for a specific data center. Since CPUs and GPUs are expensive and expected to run continuously in many data centers, the data center planners assume the computation load is fixed. As a result, cooling becomes one area where grid operators can execute demand response programs, making it the primary focus of our model. This approach allows demand response to be effectively encapsulated within the energy management system for the data center. For data center investors, the model is useful during the planning phase. It not only supports direct decisions such as the number of chillers needed for the mechanical refrigeration system, but also enables other important investment decisions such as whether to invest in additional insulation or install backup storage systems.

Traditional methods for estimating operational plan of a data center cooling have seen incremental improvements, but often do not take advantage of stochastic models for electricity prices. As a result, they can lead to increased energy use, especially during high-price periods when active cooling is most expensive. While some existing control strategies rely on fixed rules to capture the effect of electricity price (e.g., more cooling during historically low-price hours), such heuristics often cannot fully capture real-time electricity price variability.

In this paper, we model electricity prices as a stochastic process that adapts to real-time fluctuations, and develop a \textit{Markov Decision Process (MDP)} framework to optimize the operational planning of data center cooling system. We begin by applying \textit{Quantile Fourier regressions (QFR)} to capture periodic patterns in electricity prices, developing a metric that classifies them into distinct price regimes. These regimes define the states of our MDP, allowing us to generate adaptive action plans for how to run the cooling system in a data center based on real-time price realizations.

Our model includes:

\begin{itemize}

  \item \textbf{Realistic thermal dynamics} of indoor temperature changes, capturing contributions from computational heat load, passive exchange with the outdoor environment, and the cooling decisions.

  \item \textbf{Electricity price regimes}, represented as a time-inhomogeneous stochastic process using Quantile Fourier Regression (QFR) fitted to historical data of electricity prices, allowing the model to adapt to real-time prices.

  \item \textbf{Physical specifications} of the data center, including geometry and material properties, to accurately estimate thermal inertia.

\end{itemize}


We formulate the optimization problem of finding the operational plan of cooling as a linear program \cite{bertsekas2012dp}. The objective is to minimize total expected cost, including energy use and penalties for temperature violations.

Through simulations using real data for outdoor temperature and electricity prices, we show that our model yields cooling operational plans that are both cost-effective and thermally safe, outperforming conventional operational approaches.

The primary contributions of this paper include:

\begin{itemize}
    \item Using Quantile Fourier Regression (QFR) with an MDP-based cooling control framework for the first time, allowing adaptive decision-making based on real-time electricity price regimes. In our model, the data center's cooling decisions respond dynamically to price, making it a more realistic approach to optimizing market efficiency.
    
    \item Developing a model that estimates a data center's operational cooling plan, providing insights for both grid operators and data center investors; helping demand response programs and investment-phase decisions.

    \item Presenting a general model that can be tailored to geographical and physical specifications of any data center.
    
\end{itemize}


The remainder of this paper is organized as follows. Section~\ref{sec:model} explains the model formulation, Section~\ref{sec:setups} illustrates the model setups, Section~\ref{sec:results} shows the results of simulations, comparing our proposed method with conventional approaches, and Section~\ref{sec:conclu} summarizes our findings and provides suggestions for future work.

\section{Model Formulation}

\label{sec:model}
\vspace{-0.2em}
We model the data center cooling optimization problem as a discrete-time, finite-horizon Markov Decision Process (MDP). The key stochastic component of our model is the electricity price. The idea is simple: as long as the indoor temperature of the data center remains within safe limits, we aim to schedule cooling during periods of low electricity prices. To do this, we must first ask a basic question: What qualifies as a low price? This is not straightforward, since a price considered low on a summer day might be categorized as high in winter. This comes back to the fact that electricity prices exhibit stochastic and time-dependent behavior.

Analysis of historical electricity price data demonstrates a stochastic process with clear daily and seasonal periodicity. Building upon our prior work in \cite{QFR}, we use Quantile Fourier Regression (QFR) to categorize electricity prices into different quantiles and establish a clear, time-aware measure for what we mean by a high or low price. QFR is a technique that captures seasonal and diurnal periodic patterns in time series data by providing smooth, time-dependent quantile fits. We use it here to classify electricity prices in a seasonally-aware manner. Any realization of the real-time electricity price will fall into one of the discretized quantile regimes, e.g., most expensive 5\% for the season/time of day. The resulting price regimes become a part of our MDP framework. We then propose a prescribed operational plan for each possible realization in advance, derived from solving an optimization problem. For a detailed explanation of the model, see \cite{QFR}.

\begin{figure}[!t]
\centering
\includegraphics[width=0.48\textwidth]{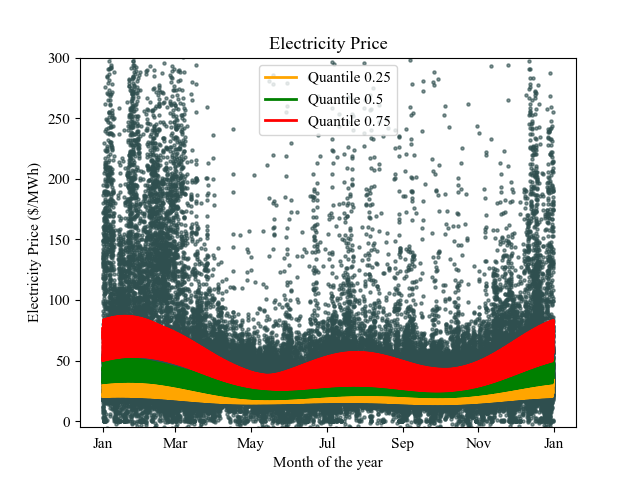}
\vspace{-0.6em}
\includegraphics[width=0.48\textwidth]{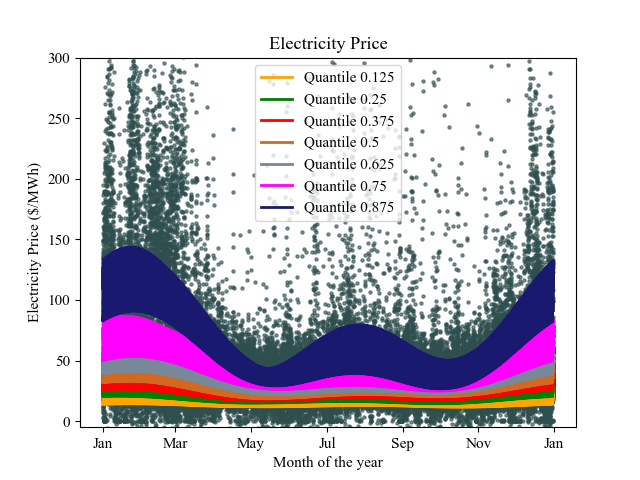}
\vspace{-0.5em}
\caption{Superimposed Quantile Fourier Regression (QFR) fits for electricity price in New England, 2011--2020, using both seasonal and daily regressors. Upper panel: Four-regime. Lower panel: Eight-regime. Both used for our QFR-MDP model.}
\label{fig:qfr_price_fits}
\end{figure}

Fig.~\ref{fig:qfr_price_fits} illustrates two examples of QFR fits to ISO New England electricity price data from 2011 to 2020, superimposed on each other (i.e. phase folded). The upper panel shows quartiles, while the lower has finer quantiles with eight regimes. The final results of our model will provide an estimate of the cooling costs and operational plan for the data center. Next, we define the components of our MDP model.

\subsection{Time horizon}

We consider an infinite time horizon MDP, which obtains the optimal policies to undertake in a stationary state (cyclostationary).
The time is $t=1,2,\ldots,N$ and wraps up after $t=N$ to $t=1$ to complete the cycle.
    
\subsection{Decision time points}

The decision time points are the same as time $t$ (e.g., every hour).

\subsection{States}

The system state at time $t$ is the tuple $s=(\theta, p)$.

\subsection{Actions}

At each time step, in state $s$, the controller determines the cooling level by selecting the number of active chillers $a$; $a=0$ denotes no chillers are active.

\subsection{State Transitions}
The state transition has two separate components that work independently: ${(\theta \rightarrow \theta',p \rightarrow p')}$
\begin{itemize}
    \item \textbf{Indoor temperature transition $\theta \to \theta'$:}  
    This is modeled as a deterministic process using an exponential heat transfer model. The indoor temperature evolves according to the following equation and is rounded to the nearest element in $\Theta$:
    \begin{equation}
    \begin{split}
    \theta_{t+1} =\; & T^{\text{out}}_t + \frac{Q_t - \eta a_t}{\gamma_{\text{env}}} \\
    & + \left( \theta_t - T^{\text{out}}_t - \frac{Q_t - \eta a_t}{\gamma_{\text{env}}} \right)
    \cdot \exp\left( -\frac{\gamma_{\text{env}} \cdot \Delta t}{C_{\text{heat}}} \right)
    \end{split}
    \end{equation}
    
    The next indoor temperature $\theta_{t+1}$ depends on:
    \begin{itemize}
        \item current temperature $\theta_t$,
        \item heat load from IT workload $Q_t$,
        \item cooling effect $\eta a_t$, representing temperature reduction resulting from the chosen cooling level $a_t$
        \item passive thermal exchange with the external environment with ambient temperature $T^{\text{out}}_t$.
        \item and building thermal properties $(C_{\text{heat}}, \gamma_{\text{env}})$.
    \end{itemize}
    
    The heat load from IT workload at time $t$ is modeled as:
    \[
    Q_t = Q^{\text{base}} + \phi \cdot i_t
    \]

    And the total heat capacity $C_{\text{heat}}$ is modeled as:
    \[
    C_{\text{heat}} = C_{\text{air}} + C_{\text{slab}} + C_{\text{equipment}}
    \]
    where:
    \begin{itemize}
        \item $C_{\text{air}} = \rho_{\text{air}} \cdot c_{p,\text{air}} \cdot V_{\text{room}}$
        \item $C_{\text{slab}} = A_{\text{floor}} \cdot h_{\text{slab}} \cdot \rho_{\text{concrete}} \cdot c_{p,\text{concrete}}$
    \end{itemize}
    These components are fixed for a given data center layout and material properties.

    \item \textbf{Electricity price regimes transition ${p \rightarrow p'}$}:
    This transition is assumed to be stochastic.
    The price regimes evolve independently via a time-inhomogeneous Markov chain with a known transition matrix $\mathbb{P}^{t}_{p \rightarrow p'}$, derived by solving a maximum likelihood estimation problem using historical data and Quantile Fourier Regression (QFR) fits. For a detailed explanation, see \cite{QFR}.
\end{itemize}


Hence, the full state transition probability from $s = (\theta, p)$ to $s' = (\theta', p')$ under action $a$ at time $t$ is:
\[
P_{t, s, s', a} =
\begin{cases}
\mathbb{P}^{t}_{p \rightarrow p'} & \text{if } \theta' = \text{round}(f(t,\theta, a)) \\
0 & \text{otherwise}
\end{cases}
\]
where $f(\cdot)$ is defined by the exponential thermal dynamics equation above.

\subsection{Cost Function}

The indoor temperature must remain within safe limits, with allowable bounds of $T_{\min}$ and $T_{\max}$ . The total cost function includes two components: the energy cost for cooling and the penalties incurred when the indoor temperature exceeds the safe range. The cost incurred for each state-action pair $(t, s, a)$ is given by:

\begin{equation}
c_{t,s,a} =\; e_{t,a} \cdot \left( \frac{\pi_{t,p}}{1000} \right)
 + \lambda_{\text{over}} \cdot v^{\text{over}}_{t+1} 
+ \lambda_{\text{under}} \cdot v^{\text{under}}_{t+1},
\end{equation}

\subsection{Optimization Problem}
We now introduce the linear program \eqref{eq:MDPLP} to compute the optimal cooling policy for the data center.

\begin{equation}
\begin{array}{lll}
\displaystyle \min_{x} & \frac{1}{N} \sum\limits_{t \in \mathcal{T}} \sum\limits_{s \in \mathcal{S}} \sum\limits_{a \in \mathcal{A}} c_{t,s,a} \cdot x_{t,s,a} & \\[10pt]
\text{s.t.} 
& \sum\limits_{s \in \mathcal{S}} \sum\limits_{a \in \mathcal{A}} x_{t s a} = 1 & \forall t \in \mathcal{T} \\[10pt]
& \sum\limits_{a \in \mathcal{A}} x_{(t+1),s',a} 
= \sum\limits_{s \in \mathcal{S}} \sum\limits_{a \in \mathcal{A}} P_{t,s,s',a} \cdot x_{t,s,a} & \forall t, s' \\[10pt]
& x_{t s a} \geq 0 & \forall t, s, a
\end{array}
\label{eq:MDPLP}
\end{equation}

where the decision variables $x_{t,s,a}$ represent the probability of being in state $s=(\theta, p)$ at time $t$ and taking action $a$.
The components of the optimization problem are as follows in order:

\begin{itemize}
    \item \textbf{Objective:} Minimizes the expected average cost of cooling the data center over the planning horizon, where the cost parameter $c_{t,s,a}$ includes electricity cost and temperature violation penalties.
    
    \item \textbf{Normalization constraint:} Ensures that at every time step $t$, sum of the probabilities is equal to 1.
    
    \item \textbf{Flow constraint:} Guarantees that the probability of reaching any future state $s'$ at time $t+1$ is consistent with the transition probabilities from all prior state-action pairs at time $t$.
    
    \item \textbf{Non-negativity:} Makes sure that all decision variables will be non-negative probabilities.
\end{itemize}

\section{Case Setups \& Policy Analysis}
\label{sec:setups}
\vspace{-0.2em}
We evaluate the proposed MDP-based cooling optimization framework using realistic traces for outdoor temperature, electricity prices, and server utilization.

\subsection{Model Setup}

We run a model of a data center cooling environment over a summer with hourly resolution ($N = 2208$). The following specifications are used:

\begin{itemize}

    \item \textbf{Data center layout:} The facility has a floor area of $3000$ m\textsuperscript{2} and ceiling height of $4$ m.

    \item \textbf{Thermal capacitance:} The total heat capacity is modeled as
    \[
    C_{\text{heat}} = C_{\text{air}} + C_{\text{slab}} + C_{\text{equipment}},
    \]
    where the air, concrete slab, and IT equipment components are computed from standard ASHRAE values.

    \item \textbf{Heat load:} Composed of a base load of $Q^{\text{base}}=1$ MW and a dynamic number of active CPU cores $i_t$ which follows a synthetic dataset for a data center, with $\phi = 10$ W/core.

    \item \textbf{Outdoor temperature:} We use 2024 hourly historical temperature data from a weather station in Boston, MA, USA.
    
    \textbf{Electricity price regimes ($\pi_{t,p}$):}  
    We fit two separate 4-state and 8-state time-inhomogeneous Markov chains to historical New England electricity prices using QFR. Each state corresponds to a price regime (e.g., low, medium, high, spike).

    \item \textbf{Cooling efficiency:} The coefficient of performance (COP) varies linearly with outdoor temperature, decreasing from $5.0$ at $15^\circ$C to $2.5$ at $40^\circ$C as follows:
        \[
    \text{COP}_t =
    \begin{cases}
    5.0 & \text{if } T^{\text{out}}_t \leq 15 \\
    2.5 & \text{if } T^{\text{out}}_t \geq 40 \\
    5.0 - \frac{T^{\text{out}}_t - 15}{10} & \text{otherwise}
    \end{cases}
    \]
    
\end{itemize}

\subsection{Case Setup}
For the cooling system, we assume having 4 chillers with $1.25 MW$ power. We compare two different price regime models as follows. This setup is designed to assess whether increasing the number of electricity price regimes provides meaningful improvements in cost-aware cooling policies or only increases model complexity.

\textbf{Price regime setup:}
\begin{itemize}
    \item A 4-regime model where each price regime represents a coarse level of electricity cost (e.g., low, medium, high, spike).
    
    \item An 8-regime model with finer granularity, capturing more detailed variations in electricity prices.
\end{itemize}

\subsection{Policy Visualization}

\begin{figure}[!t]
\centering
\includegraphics[width=0.48\textwidth]{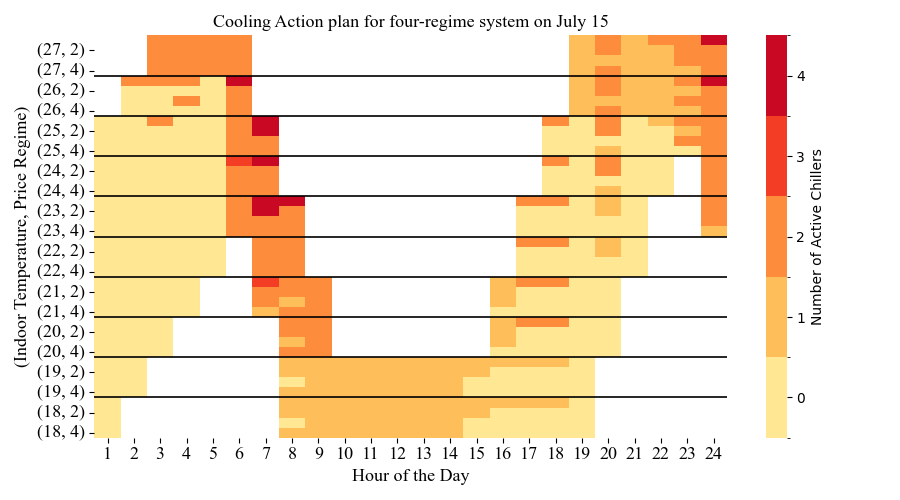}

\includegraphics[width=0.48\textwidth]{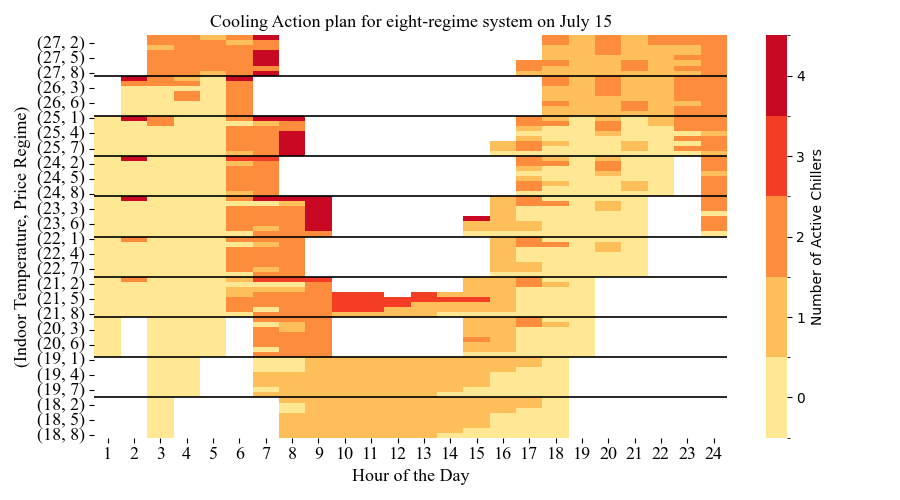}
\vspace{-1.2em}
\caption{Optimal daily cooling action plan under the QFR-MDP electricity-price model for a data center with four independent chillers (example day: July 15). Colors indicate the number of active chillers: \textbf{red}—all four chillers active, \textbf{orange}—two chillers active, \textbf{yellow}—no chillers active. Upper panel: Four-regime. Lower panel: Eight-regime.}
\label{fig:cooling_plan_4chillers}
\end{figure}

By solving \eqref{eq:MDPLP}, we obtain the optimal action for each hour of the day, tailored to each specific day of the summer. Fig.~\ref{fig:cooling_plan_4chillers} illustrates the optimal cooling action-plan $a_t$ as the number of active chillers selected at each time step $t$ under different electricity price regimes and indoor temperatures for a specific example day in the summer (July 15th). The policy exhibits anticipatory behavior: when high-price regimes are expected, the system pre-cools during cheaper hours to avoid expensive cooling during peak prices. Two key points should be considered. First, the system prefers to remain at higher temperatures as long as they stay within safe limits, since there is no additional benefit in cooling beyond safe limits. Second, pre-cooling cannot push the temperature below the minimum allowed limit. Both of these points are clearly observable in Fig.~\ref{fig:cooling_plan_4chillers}.

\section{Numerical Results}
\label{sec:results}
\vspace{-0.2em}

In this section, we validate the performance of our proposed QFR-based MDP model through simulations, using real-time electricity price data from 2011 to 2024 from ISO New England, and outside temperature data from the Boston, MA, station. The simulation is conducted in-sample from 2011 to 2020 and out-of-sample from 2021 to 2024. In order to do that, we introduce two well-known heuristic methods of cooling as baselines.

We compare the proposed QFR-MDP policy with two benchmarks:
\begin{itemize}
  \item \textbf{Greedy Controller:} Selects the minimum cooling action such that the next temperature does not exceed \( T_{\max}= \)27°C.
  \item \textbf{Peak-Hour Based Fixed Rule:} No cooling during traditional peak hours (4 PM to 7 PM), without adapting to real-time price.
\end{itemize}

\begin{figure}[!t]
\centering
\includegraphics[width=0.48\textwidth]{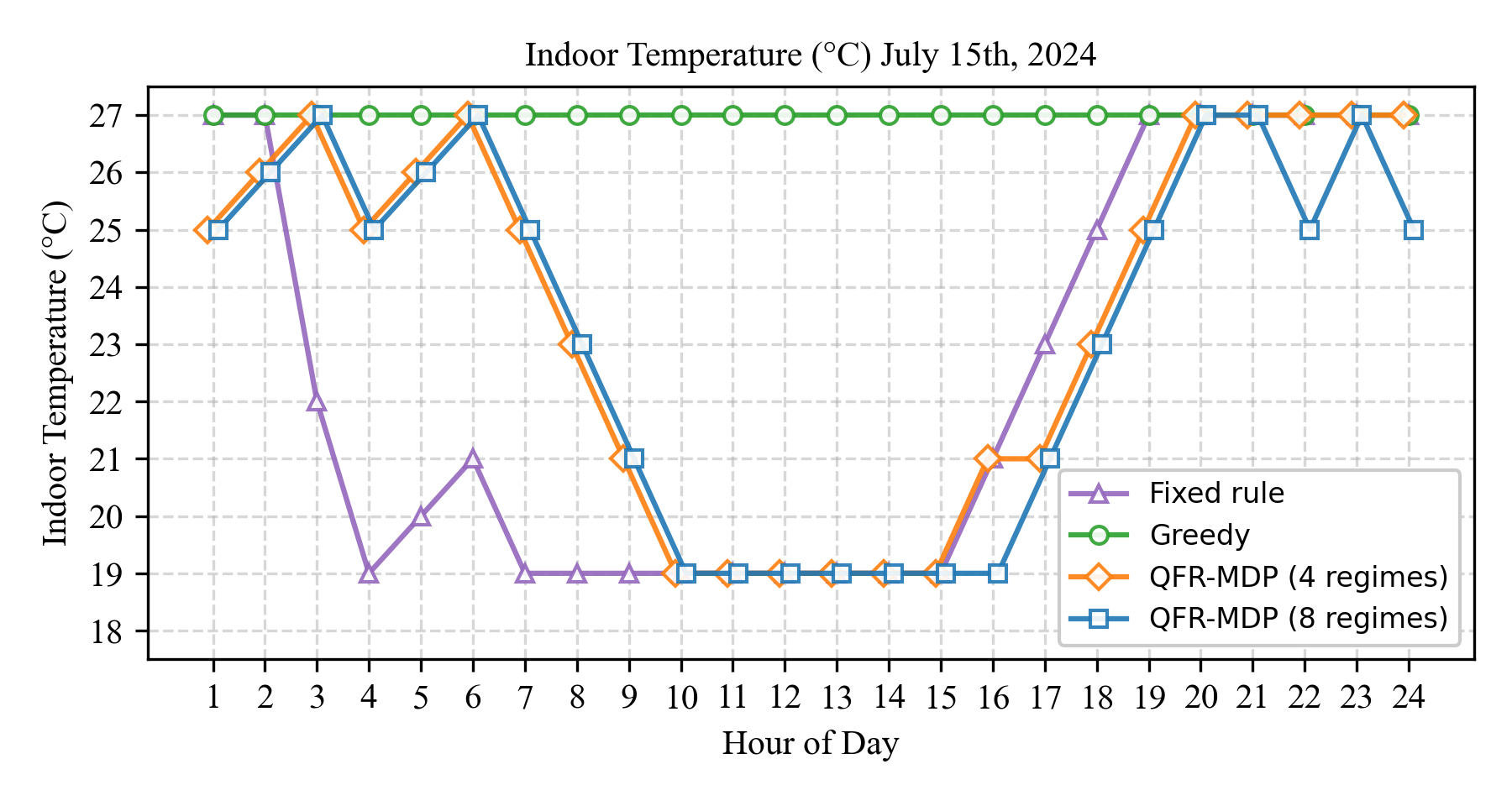}
\vspace{-1.6em}
\caption{Inside temperature across 24 hours on July 15, 2024, under three control methods: Fixed rule, Greedy, and QFR-MDP. All three methods maintain the indoor temperature within the safe limits, 18–27°C.}
\label{fig:indoor_temperature}
\end{figure}

\begin{figure}[!t]
\centering
\includegraphics[width=0.48\textwidth]{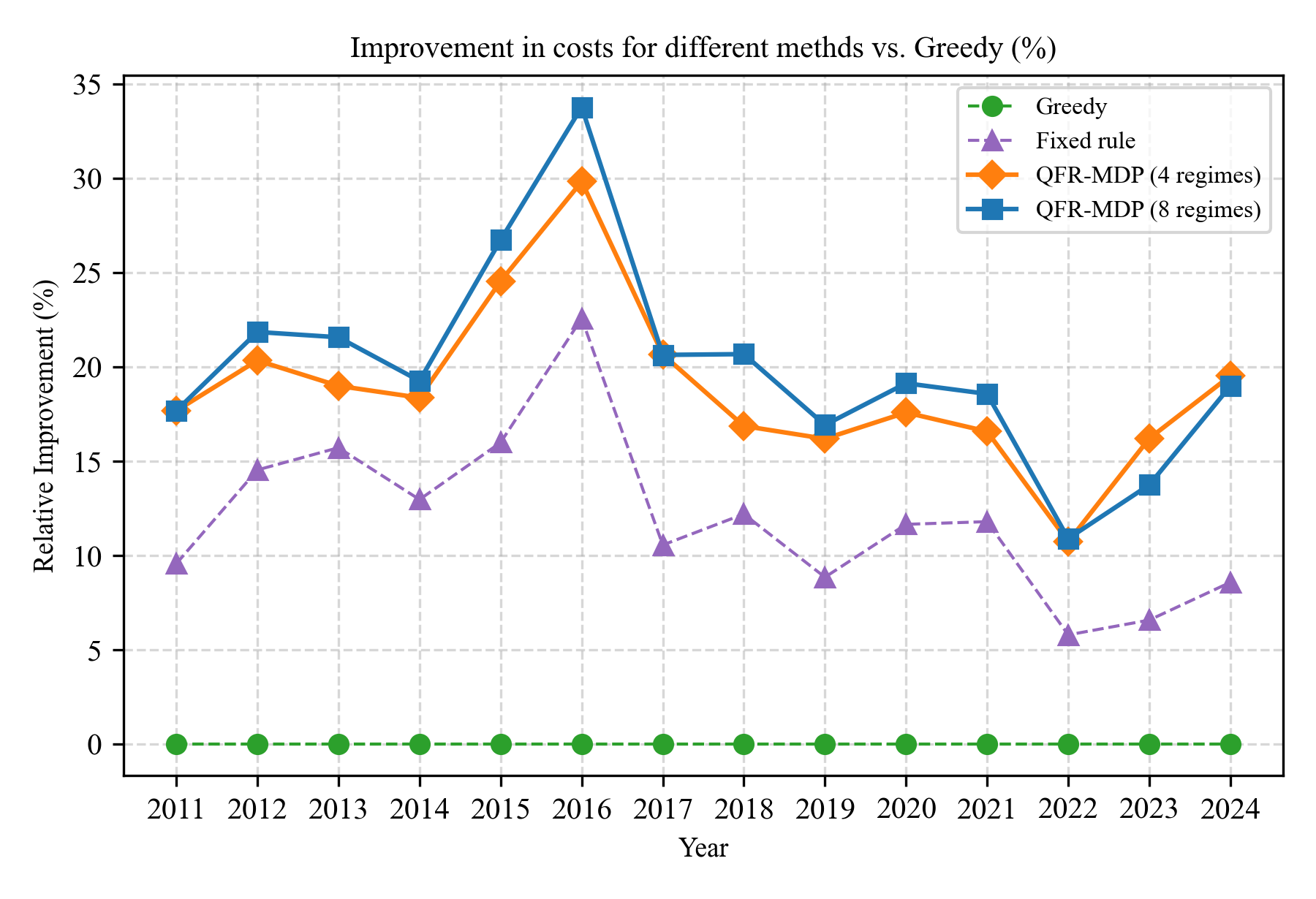}
\vspace{-1.6em}
\caption{Relative improvement in cooling cost across control methods from 2011 to 2024 summers. The QFR-MDP approach consistently achieves lower cooling costs compared to the Greedy as a baseline.}
\label{fig:costs}
\end{figure}

Fig.~\ref{fig:indoor_temperature} shows the indoor temperature on a sample date (July 15th, 2024) under different policies. In all cases, the system maintains temperature within safe bounds, 18–27°C. The Greedy method holds the temperature at the upper limit \( T_{\max}= \)27°C. The Fixed rule aggressively cools in the middle of the night (2 AM - 3 AM) and turns off during peak hours (3 PM - 7 PM). In contrast, our proposed QFR-MDP method dynamically adapts cooling to stay within the comfort range based on real-time prices.
Fig.~\ref{fig:costs} compares the energy cost of cooling the data center under different policies in a simulation across multiple years of summers. The QFR-MDP method consistently achieves lower energy costs compared to both baselines. This stems from the inability of the other two methods to react to stochastic variability in electricity prices. These results demonstrate the potential of QFR-MDP scheduling for electricity price-aware, thermally safe data center operations. The fact that the eight-state system achieves lower costs than the four-state system in most summers suggests that using finer bins can lead to more effective policies in terms of average cost.

\section{Conclusion}
\label{sec:conclu}
\vspace{-0.2em}

We present an MDP-based model for electricity price-aware cooling in a data center, utilizing QFR fits to categorize real-time prices into different price regimes. These fits are time-dependent, enabling our model to capture both the diurnal and seasonal patterns in the electricity price. This model is specifically designed for investment phase decisions by providing an estimate of the operational costs and plan for cooling a data center. It offers an accurate estimate of optimized cooling costs well before the construction stage of a data center, which is important for both engineering design and financial planning.

Over 14 years of summer simulations using real-time electricity price data, we show that our model consistently has lower cooling costs compared to other heuristic baselines.
Future work may explore the co-optimization of price-aware cooling and workload scheduling. Another direction is to extend our framework to incorporate broader sustainability objectives, such as reducing carbon emissions from electricity use while considering fairness (see \cite{sahebdel2025lead,sahebdel2024holistic}).

\section*{Acknowledgment}
\vspace{-0.2em}
This material is based upon work supported by the U.S. Department of Energy’s Office of Energy Efficiency and Renewable Energy (EERE) under the Wind Energy Technologies Office (WETO) Award Number DE-EE0011269, the Massachusetts Clean Energy Center and the Maryland Energy Administration. The views expressed herein do not necessarily represent the views of the U.S. Department of Energy, the United States Government, the Massachusetts Clean Energy Center or the Maryland Energy Administration.

\bibliographystyle{IEEEtran}
\bibliography{References}

\end{document}